\documentclass[conference]{IEEEtran}
\IEEEoverridecommandlockouts
\usepackage{cite}
\usepackage[utf8]{inputenc} 
\usepackage[T1]{fontenc}    
\usepackage{amsmath,amssymb,amsfonts}
\usepackage{booktabs}
\usepackage{multirow}
\usepackage{algorithmic}
\usepackage{graphicx}
\usepackage{textcomp}
\usepackage{xcolor}
\usepackage{nomencl}
\usepackage{comment}
\def\BibTeX{{\rm B\kern-.05em{\sc i\kern-.025em b}\kern-.08em
    T\kern-.1667em\lower.7ex\hbox{E}\kern-.125emX}}
\begin{document}

\title{Optimal Battery Charge Scheduling For Revenue Stacking Under Operational Constraints Via Energy Arbitrage}

\author{\IEEEauthorblockN{Alban Puech$^{1,3}$, Gorazd Dimitrov$^1$, Claudia D'Ambrosio$^2$} 
\IEEEauthorblockA{\textit{$^1$ Institut Polytechnique de Paris, Palaiseau, France}\\
\textit{$^2$ LIX CNRS, École polytechnique, Palaiseau, France}\\
\textit{$^3$ DEIF Wind Power Technology, Klagenfurt, Austria}\\
\{alban.puech, gorazd.dimitrov, claudia.dambrosio\}@polytechnique.edu}}


\maketitle

\begin{abstract}
    As the share of variable renewable energy sources increases in the electricity mix, new solutions are needed to build a flexible and reliable grid. Energy arbitrage with battery storage systems supports renewable energy integration into the grid by shifting demand and increasing the overall utilization of power production systems. In this paper, we propose a mixed integer linear programming model for energy arbitrage on the day-ahead market, that takes into account operational and availability constraints of asset owners willing to get an additional revenue stream from their storage asset. This approach optimally schedules the charge and discharge operations associated with the most profitable trading strategy, and achieves between 80\% and 90\% of the maximum obtainable profits considering one-year time horizons using the prices of electricity in multiple European countries including Germany, France, Italy, Denmark, and Spain.
\end{abstract}

\section{Introduction}
In IEA's pathway to net-zero CO2 emissions by 2050 \cite{iea_2021}, electricity accounts for 50\% of energy use in 2050, with 70\% produced by wind and solar photovoltaic. As reliance on variable energy sources (VES) increases, important investments in grid-scale storage solutions are required to balance their variable and non-dispatchable aspects, and to increase grid flexibility, capacity, and reliability \cite{pavarini_2018}. While Pumped-storage hydroelectricity (PHS) represents 90\% of the world's storage \cite{MCWILLIAMS2022147}, battery system storage (BSS) are expected to see the largest market growth and can play a similar role \cite{schoenfisch_dasgupta_2022}. BSS have many use cases, from replacing backup emergency power units in hospitals to frequency response systems. Some of the applications of BSS are listed in \cite{9509501}, including energy arbitrage, which helps accommodate changes in demand and production by charging BSS during low-cost periods and discharging them during consumption peaks when prices are high. This assists the integration of renewable energy and allows higher utilization of VES by shifting demand to align with important production periods. The profitability of energy arbitrage depends on factors such as battery cost and price differences between low and high-demand periods. Intraday price variability represents an opportunity for asset owners to stack additional revenue streams obtained from arbitrage to their existing ones, providing additional incentives for operating and investing in grid-scale BSS. Furthermore, Energy arbitrage increases the return on investment (ROI) of BSS. Paper \cite{karaduman2020economics} shows other positive societal impacts, such as increasing the return to renewable production and reducing CO2 emissions. Energy arbitrage with BSS requires a charging/ discharging schedule. This schedule needs to consider battery characteristics, such as the charging curve, maximum cycles per day, and discharging rate. Many papers have dealt with this problem using model predictive control based on linear programming (LP) or mixed-integer linear programming (MILP) optimizations \cite{Krishnamurthy2018EnergySA,8340703,6477197,6168818,7563833}. However, those approaches often assume constant charging rates and no battery degradation.

In this paper, we propose an approach to generate optimal charging schedules for arbitrage on the Day Ahead Market (DAM) under availability constraints that considers variable charging rates and battery degradation. Our contributions are:
\begin{itemize}
\item A battery model with variable charging rates, discharge efficiency decrease, and capacity fading, that can be parametrized by the user.
\item A MILP formulation of the profit maximization task via energy arbitrage on the DAM.
\item A Python library that outputs daily optimal schedules from user inputs using built-in price forecasting, and allows running the algorithm on historical prices.
\item A quantitative comparison of the performances of our predictive optimization model relying on price forecasts with a baseline that outputs optimal schedules a posteriori, for different key hyperparameters and different European countries.
\item A discussion of the impact of capacity fading, discharge efficiency decrease, and charging rates on profits.
\end{itemize}
Our battery model and optimization framework are introduced in Section \ref{methodology}. The main results: the comparison of our output schedules to the optimal a posteriori schedule, the performances of our algorithm on datasets considering a one-year time horizon in multiple European countries, and the discussion of the impact of prediction accuracy and battery degradation on the obtained profits are presented in Section \ref{results}. In Section \ref{conclusion} we conclude and discuss future work. 


\begin{figure}[!htbp]
\centering
\includegraphics[width=\columnwidth]{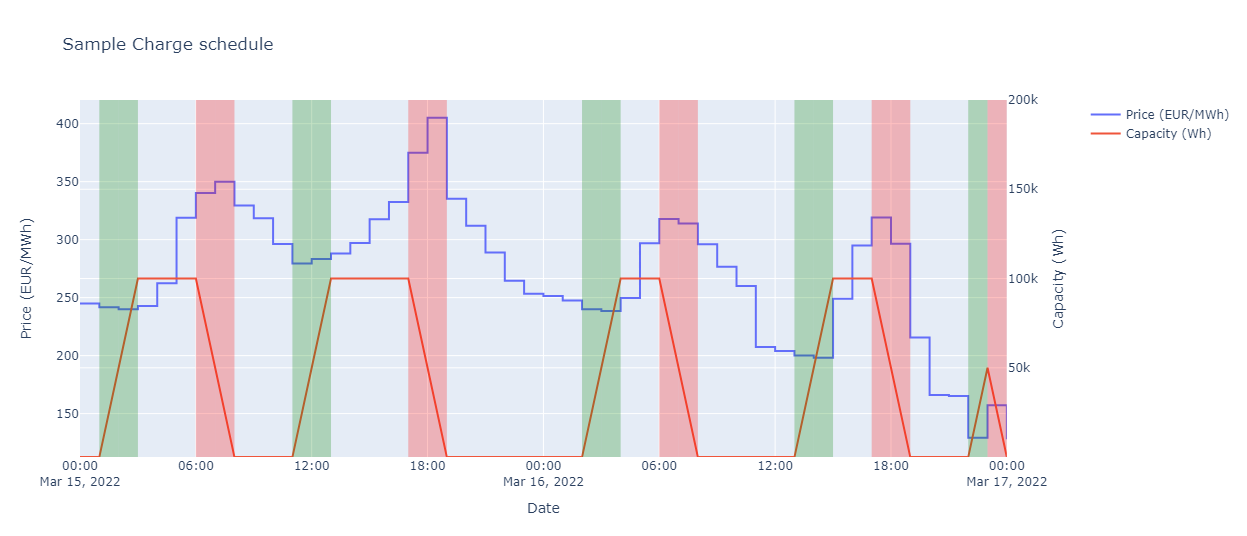}
\caption{\label{fig:schedule}Sample schedule for a 100 kWh battery, assuming no grid cost, a constant charge rate of 0.5W/Wh, and without any limit on the number of charging/discharging cycles}
\end{figure}

\nomenclature{$Q_{0}$}{Initial energy capacity (before capacity fading) in Wh}
\nomenclature{$\text{SOC}(h)$}{State of charge at time $h$}
\nomenclature{$\Delta E(h)$}{Change in energy stored in the battery between hour $h$ and hour $h+1$ in Wh
}
\nomenclature{$\overline{\Delta}(\text{SOC}(h))$}{Maximum positive energy change between hour $h$ and hour $h+1$ in Wh}
\nomenclature{$\underline{\Delta}((\text{SOC}(h)))$}{Minimum negative energy change between hour $h$ and hour $h+1$ in Wh}
\nomenclature{$\Delta E_{day}(d)$}{Total amount of energy exchanged during day $d$}
\nomenclature{$f_c(soc)$}{Charging rate at $\text{SOC}=soc$ in W/Wh}
\nomenclature{$f_d(soc)$}{Discharging rate at $\text{SOC}=soc$ in W/Wh}
\nomenclature{$Q(d)$}{Energy capacity on day d in Wh}
\nomenclature{$n_{cycles}(d)$}{Number of battery cycles done up to day $d$}
\nomenclature{$\alpha$}{Battery capacity decay}
\nomenclature{$\eta(d)$}{Battery charge efficiency on day $d$}
\nomenclature{$\beta$}{Battery charge efficiency decay}
\nomenclature{$p_d(h)$}{Electricity price at hour $h$ on day $d$ in EUR/Wh}
\nomenclature{$\Tilde{p}_d(h)$}{Prediction of the electricity price at hour $h$ on day $d$ in EUR/Wh}
\nomenclature{$l$}{Numbers of days used for the computation of the forecasted price}

\printnomenclature

\section{Methodology} 
\label{methodology}
\subsection{Desired algorithm output}

In a DAM, the prices for the next day are determined through an auction process that typically closes before midday. The auction sets prices simultaneously for all 24 hours of the next day, allowing for the acquisition of buy-sell pairs without leaving any open positions. This means that our algorithm has to generate the schedule $S_d$ for the day $d$ before market closure on day $d-1$. In Fig. \ref{fig:schedule}, a sample schedule is shown. We denote by $E_d(h)$ the change in energy stored in the battery between hour $h$ and hour $h+1$ on day $d$. Our scheduling problem is, thus, formulated as an optimization task that outputs a sequence $S_d$ of length 24 (one set point per hour): $$S_d = (E_d(0), E_d(1), \dots, E_d(23))$$

 \subsection{Price prediction}
The exact prices for electricity on the day $d$ are not known until after the day-ahead auction has closed. Thus, there is a need to rely on price predictions to be able to get our schedule before the market closure. In the present approach, the mean hourly prices over the last $l$ preceding day $d-1$ are used as a prediction of the electricity prices for the day $d$. We, then, generate the schedule for the day $d$ before market closure on day $d-1$.  $\Tilde{p_d}(h)$, the forecast of the electricity price at hour $h$ denoted by $p_d(h)$ is given by $ \Tilde{p_d}(h) = \sum_{d' = d-l}^{d-1} \frac{p_d(h)}{l}$.

\subsection{Grid costs}

There are variable and fixed costs $\mathrm{vgc}_d(h)$ and  $\mathrm{fgc}_d(h)$ associated with buying and selling electricity from and to the grid that depends on the day $d$ and the hour $h$. The value of $\mathrm{vgc}_d(h)$ is expressed in EUR/MWh and is to be multiplied by the total amount of energy exchanged with the grid. The value of $\mathrm{fgc}_d(h)$ is expressed in EUR and is paid for every hour when energy has been exchanged with the grid.

\subsection{Battery model}

We introduce a battery model that simulates a battery of a given energy capacity $Q$. The state of charge at hour $h$, $\text{SOC}(h)$ is computed as stated below, where $E_{init}$ denotes the initial energy stored in the battery at hour $0$:
\begin{IEEEeqnarray*}{rCl}
\label{soc}
    \text{SOC}(h) = \begin{cases} \frac{E_{init} + \sum_{t = 0}^{h-1} E(t)}{Q} \text{ if } h \in \{1,2,..,24\}\\
    \frac{E_{init}}{Q} \text{ if } h=0.
    \end{cases}
\end{IEEEeqnarray*}

\subsection{Variable charging rates}
The state of charge $\text{SOC}$ of the battery can be computed as a function of the charging rate $r(t)$, expressed in W/Wh, which varies during the charge:
\begin{equation*}
\label{varying}
    \text{SOC}(h) = \text{SOC}(0) + \int_{t=0}^{h} r(t)dt.
\end{equation*} 

\begin{figure}
\centering
\includegraphics[width=\columnwidth]{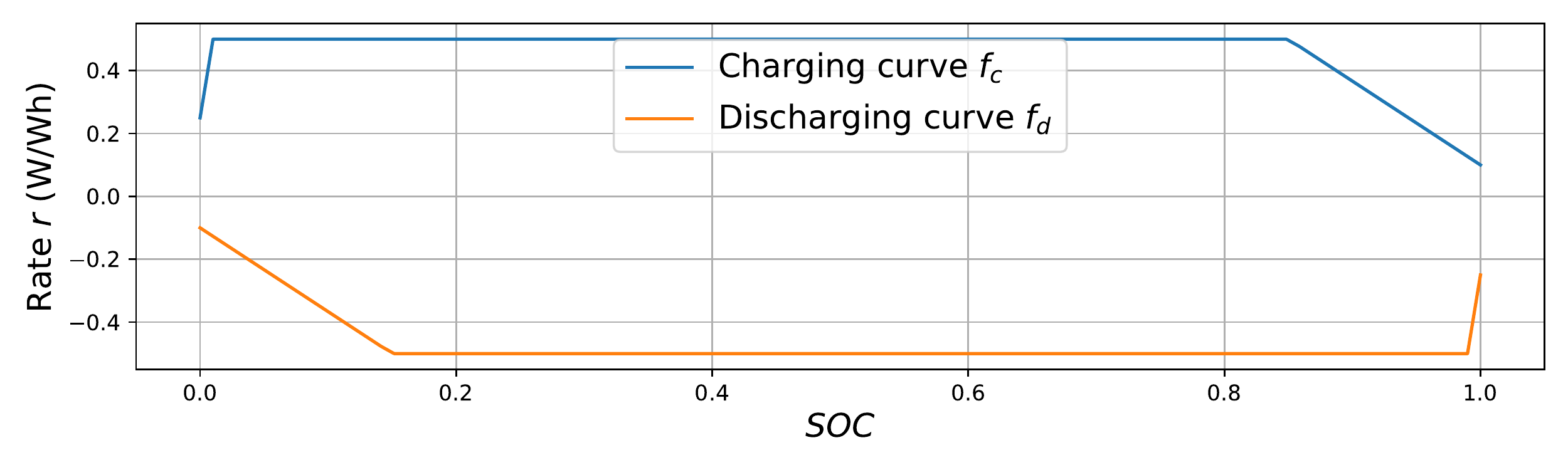}
\caption{\label{fig:charging_curve} Charging curve $f_c$ and discharging curve $f_d$ of a lithium-ion battery }
\end{figure}

The charging and discharging curves give the charging and discharging rates as a function of the $\text{SOC}$. Fig. \ref{fig:charging_curve} shows sample charging and discharging curves for a lithium-ion battery. The charging and discharging curves $f_c$ and $f_d$ can be used to derive the maximum (positive) hourly \text{SOC} change $\overline{\Delta}(\text{SOC}(h))$ and the minimum (negative) hourly \text{SOC} change $\underline{\Delta}(\text{SOC}(h))$ between hour $h$ and hour $h+1$ :
\begin{IEEEeqnarray*}{rCl}
    \overline{\Delta}(\text{SOC}(h)) = \min\bigg(&(1-\text{SOC}(h)),
    &\int_{t=h}^{t=h+1} f_c(\text{SOC}(t)) dt \bigg) \\
    \underline{\Delta}(\text{SOC}(h))  = \max\bigg(&-\text{SOC}(h),
    &\int_{t=h}^{t=h+1} f_d(\text{SOC}(t)) dt \bigg) 
\end{IEEEeqnarray*}
The constraint on the change in energy stored in the battery between hour $h$ and hour $h+1$, then, becomes:
\begin{equation*}
       Q \cdot \overline{\Delta}(\text{SOC}(h))  \leq E(h)  \leq Q \cdot \underline{\Delta} (\text{SOC}(h))
\end{equation*}

\begin{figure}
\centering
\includegraphics[width=\columnwidth]{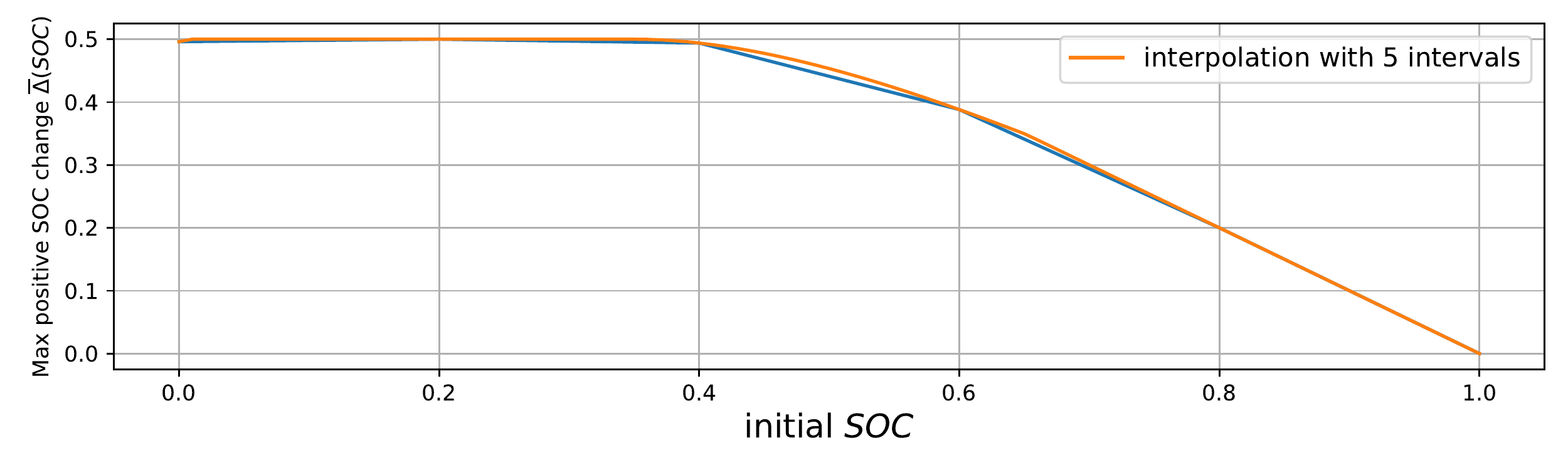}
\includegraphics[width=\columnwidth]{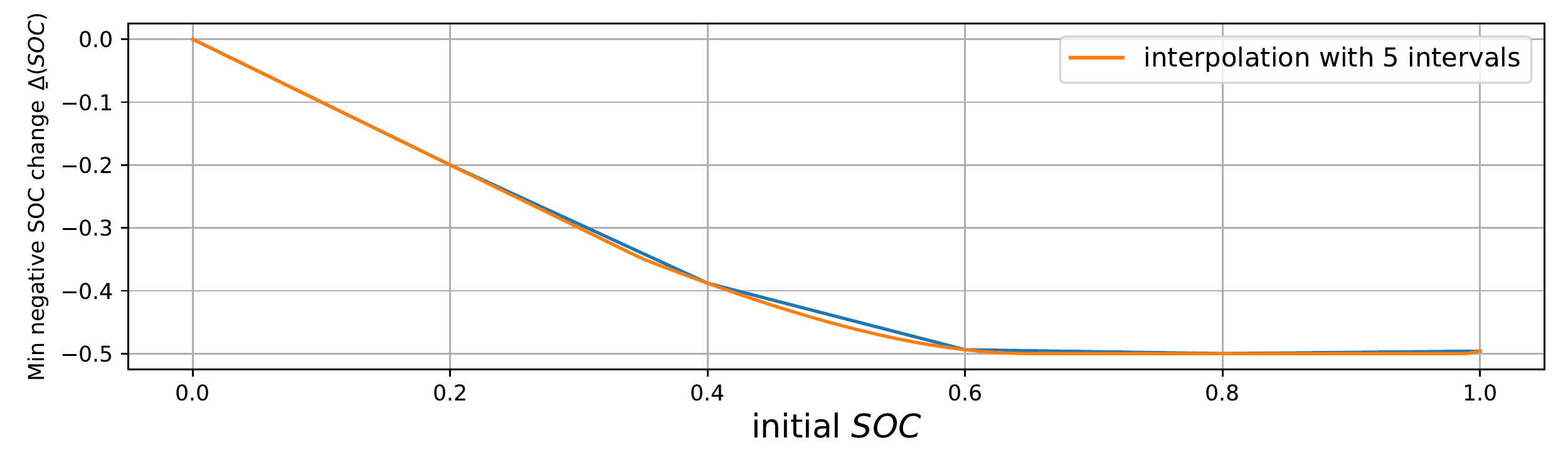}
\caption{\label{fig:energy_change} Max. positive \text{SOC} change $ \overline{\Delta}$ (top) and min. negative \text{SOC} change $ \underline{\Delta}$ (bottom) in one hour of charge as a function of $\text{SOC}$}
\end{figure}

As those functions are not linear, we use convex combination formulations to obtain piecewise linear approximations. We discretize the interval $[0,1]$ into $N_{int}$ intervals of length $N_{int}^{-1}$. This adds the following set of constraints
($\forall h \in \{0,1,\dots,23\},  \forall k \in \{1,2,\dots,N_{int}\}$):
\allowdisplaybreaks
\begin{IEEEeqnarray}{lll}
&\sum_{k=1}^{N_{int}}\bigg( \lambda^k_h \cdot \frac{k-1}{N_{int}} + \mu^k_h \cdot \frac{k}{N_{int}}\bigg) = \text{SOC}(h) \label{pwl:eq1}\\
&\lambda^k_h + \mu^k_h = y^k_h \label{pwl:eq2}\\
&\sum_{k=1}^{N_{int}} y^k_h = 1 \label{pwl:eq3}\\
& Q \cdot\sum_{k=1}^{N_{int}}\bigg(\lambda^k_h \cdot \overline{\Delta}\big(\frac{k-1}{N_{int}}\big) + \mu^k_h \cdot \overline{\Delta}\big(\frac{k}{N_{int}}\big)\bigg) \geq E(h) \label{pwl:eq4}\\
& Q \cdot \sum_{k=1}^{N_{int}}\bigg(\lambda^k_h \cdot \underline{\Delta}\big(\frac{k-1}{N_{int}}\big) + \mu^k_h \cdot \underline{\Delta}\big(\frac{k}{N_{int}}\big)\bigg)\leq E(h) \label{pwl:eq5}\\
&\mu^k_h \in [0,1], \quad \lambda^k_h \in [0,1], \quad y^k_h \in \{0,1\}  \label{pwl:eq6}
\end{IEEEeqnarray}
Constraints \eqref{pwl:eq1} are used to set the values of the convex combination weights $\lambda$ and $\mu$, \eqref{pwl:eq2} guarante the convex requirement (sum to 1), \eqref{pwl:eq3} ensure that one interval is selected for each point, and \eqref{pwl:eq4} and \eqref{pwl:eq5} set the maximum and minimum energy change for each hour. On Fig. \ref{fig:energy_change} values of $\overline{\Delta}$ and $\underline{\Delta}$ for different initial $\text{SOC}$ values are shown. The blue lines correspond to the true functions, while the orange ones correspond to a piecewise linear approximation using 5 intervals (6 cut-points).

\subsection{Availability constraints} 
 Energy arbitrage can represent an additional revenue stream for asset owners using their batteries for another main purpose (e.g. backup energy source) \cite{stacking}. In this case, the battery has to be charged enough at some key points in time, or, on the contrary, needs to be able to store energy excess, for example, from renewable electricity production systems. We define the parameters  $\overline{\text{SOC}}(h)$ and  $\underline{\text{SOC}}(h)$, which sets the upper and lower bounds for the state of charge at hour $h$. This adds the following simple bounds:
\begin{equation*}
   \underline{\text{SOC}}(h) \leq \text{SOC}(h) \leq \overline{\text{SOC}}(h) \quad \forall h \in \{1,2,\dots,24\} 
\end{equation*}

\subsection{Battery degradation}
 We characterize and model the evolution of the state of health of the battery by decreasing the initial capacity $Q$ and by decreasing the charge efficiency factor $\eta$.

\subsubsection{Capacity fading} During the life of the battery, its capacity decreases with the number of cycles \cite{fading}. This is referred to as capacity fading. We model capacity fading by decreasing its initial capacity $Q_{0}$ with a daily frequency, linearly with the number of cycles. The battery is characterized by its cycle life $Cycle_{\max}$. The battery capacity reaches 80\% of its initial capacity after $Cycle_{\max}$ cycles. The MILP parameter corresponding to the capacity on day $d$, denoted by $Q(d)$, is computed as below, where $E_{day}(d)$ denotes the total amount of energy exchanged during day $d$:
\allowdisplaybreaks
\begin{IEEEeqnarray*}{lll}
&Q(d) = max(Q_{min}, Q_{0} - \alpha \cdot n_{cycles}(d))\\
&n_{cycles}(0) = 0, \quad n_{cycles}(d) = \frac{\sum_{d' = 0}^{d-1} E_{day}(d')}{2 Q_0}, d \geq 1\\
&Q_{min} = 0.8 \cdot Q_{0} \\
&\alpha = \frac{Q_{0}-Q_{min}}{Cycle_{\max}} 
\end{IEEEeqnarray*}
\\

\subsubsection{Charge efficiency decrease} The charge efficiency $\eta$ decreases with the number of cycles. We update it similarly to the capacity:
\begin{IEEEeqnarray*}{lll}
&\eta(d) = max(\eta_{min}, \eta_0 - \beta \cdot n_{cycles}(d)   ) \\
&\eta_{min} = 0.8 \cdot \eta_0 \\
&\beta = \frac{\eta_0-\eta_{min}}{Cycle_{\max}} 
\end{IEEEeqnarray*}






\subsection{Objective function}
Our goal is to maximize forecasted profits. If the battery is charging at hour $h$, the forecasted (negative) profit $\pi(h)$ is the forecasted cost of the electricity bought added to the variable and fixed grid costs:
\begin{equation*}
    \pi(h) = -(\Tilde{p}(h) + \mathrm{vgc}(h))\cdot E(h) + \mathrm{fgc}(h).
\end{equation*}

If the battery is discharging at hour $h$, the forecasted profit $\pi(h)$ is the revenue from the electricity sold taking into account the discharge efficiency $\eta(d)$, net of the variable and grid costs:
\begin{equation*}
    \pi(h) = -(\Tilde{p}(h) + \mathrm{vgc}(h))\cdot \eta(d) \cdot E(h) + \mathrm{fgc}(h).
\end{equation*}

%
%
So as to have a linear objective function, we keep track of the chosen actions at each hour $h$. We introduce the binary variables $z(h)$ that indicate if the amount of energy stored in the battery changes during hour $h$. We divide $E(h)$ into positive $E^+(h)$ and negative parts $E^-(h)$ and we obtain the objective function $\pi_{total}$ in a linear form: 
\begin{IEEEeqnarray*}{lll}
\pi_{total} = \sum_{h=0}^{23} & & \big[(E^+(h) - E^-(h) \cdot \eta) \cdot \Tilde{p}_d(h) \\ & & + (E^+(h) + E^-(h) \cdot \eta) \cdot \mathrm{vgc}(h)) \\ 
& &+ z(h) \cdot \mathrm{fgc}(h) \big]
\end{IEEEeqnarray*}

with added constraints:
\begin{IEEEeqnarray*}{lll}
& &E^+(h) - E^-(h) = E(h) \\
& &0 \leq E^+(h) \leq M \cdot y(h), \quad 0 \leq E^-(h) \leq M \cdot (1-y(h)) \\
& &\varepsilon \cdot z(h) \leq E^+(h) + E^-(h) \leq  M \cdot z(h)\\
& &y(h), z(h) \in \{0,1\} \\
& & M = Q_0, \varepsilon = 0.01 \\
\end{IEEEeqnarray*}

\section{Results and Discussion} \label{results}
\subsection{Data and experiment setup}
We test our algorithm using the electricity prices for the year 2022 from the European Wholesale Electricity Price Dataset \cite{price}, which contains average hourly wholesale day-ahead electricity prices for European countries. We use a battery with a 1 MWh initial capacity ($ Q_{0} = 10^6$), and with charging and discharging curves as shown on Fig.\ref{fig:charging_curve}. We also set $\eta_{0} = 0.99$, $Cycles_{max}=4000$, and assume no fixed grid costs. The variable grid costs are set to 5 EUR/MWh for every day and every hour of the run, i.e. $\mathrm{vgc}(h) = 5 \cdot 10^{-6}$. We set the maximum and minimum state of charge to 0 and 1 respectively for every hour of every day of the time-horizon, except midnight, i.e. $\overline{\text{SOC}}(h) = 1, \underline{\text{SOC}}(h) = 0$. We use our method to calculate daily schedules for each day of the year 2022. Each schedule starts at 00:00 (midnight) and ends at 24:00. We set $\overline{\text{SOC}}(24) = 0 \text{ and } E_{init} = 0$ for every day so that the battery is always empty at the end of the day. We use the Amplpy library to directly integrate the AMPL \cite{AMPL} model into our Python library. We update the battery parameters (e.g. $\eta(d)$, $Q(d)$) after every call to the MILP solver (i.e. after each daily schedule is generated). We compare our algorithm described above (denoted by MILP-P) with MILP-O, where MILP-O denotes the algorithm that uses the true prices instead of price forecasts and provides us with the optimal schedule.




\subsection{Results}
Table \ref{table-metrics-algorithms} shows the results of MILP-O and MILP-P with different window sizes $l$ on the German electricity prices. The optimization for the entire year 2022 took under $1.5$ minutes for the whole year with the Gurobi solver on a Huawei Matebook with an AMD Ryzen 5 4600H and 16.0 GB RAM.

\begin{figure}[t]
\centering
\includegraphics[width=\columnwidth]{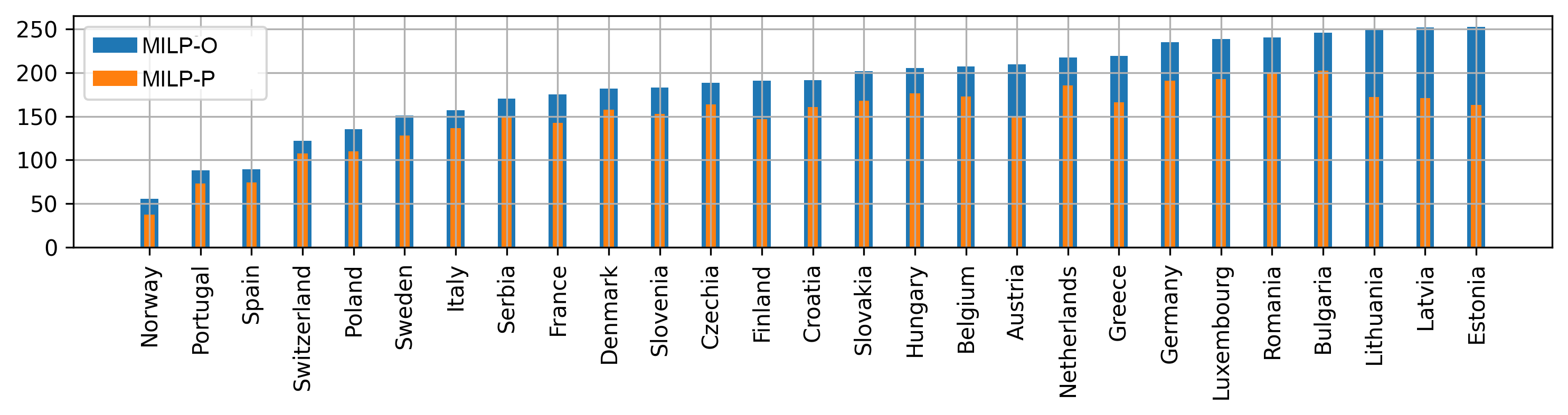}
\caption{\label{fig:profit_countries} Profit achieved on one-year-long simulations on the 2022 electricity prices of different European countries}
\end{figure}

\begin{figure}
\centering
\includegraphics[width=\columnwidth]{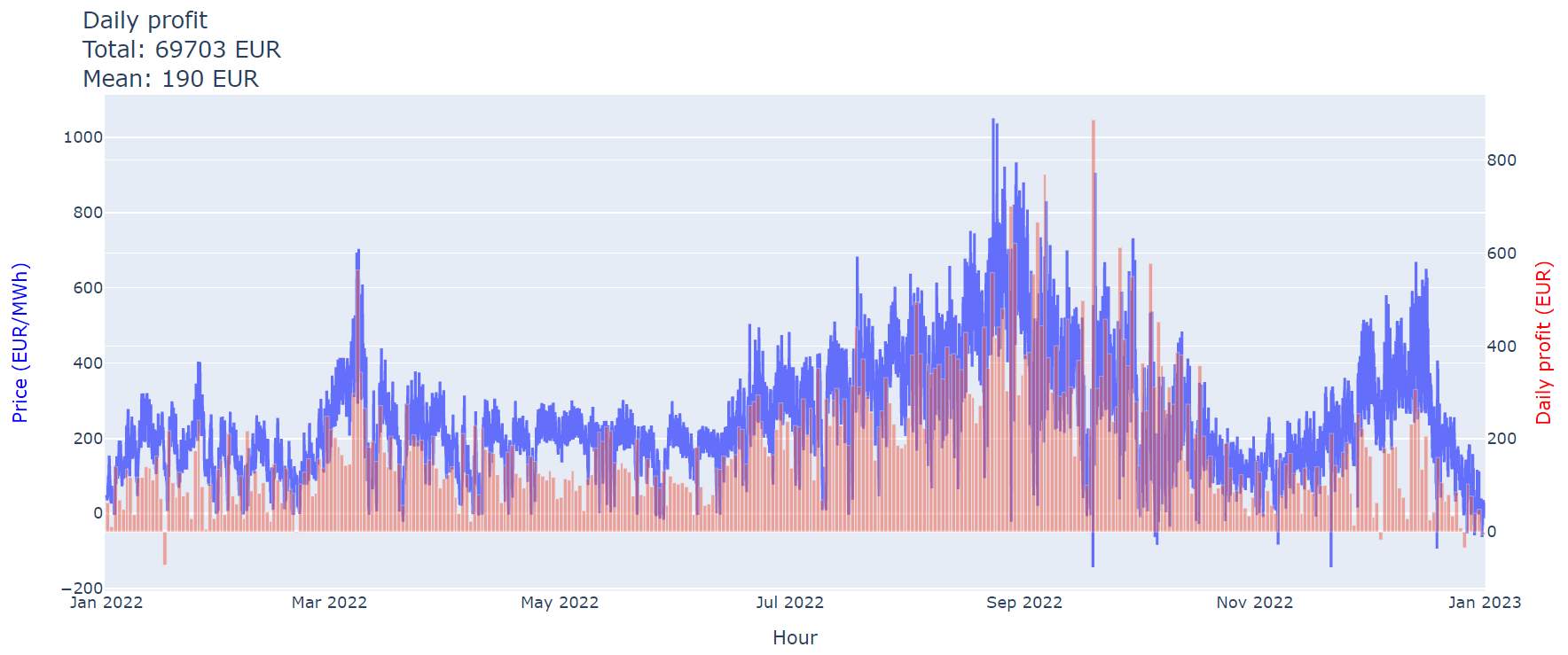}
\includegraphics[width=\columnwidth]{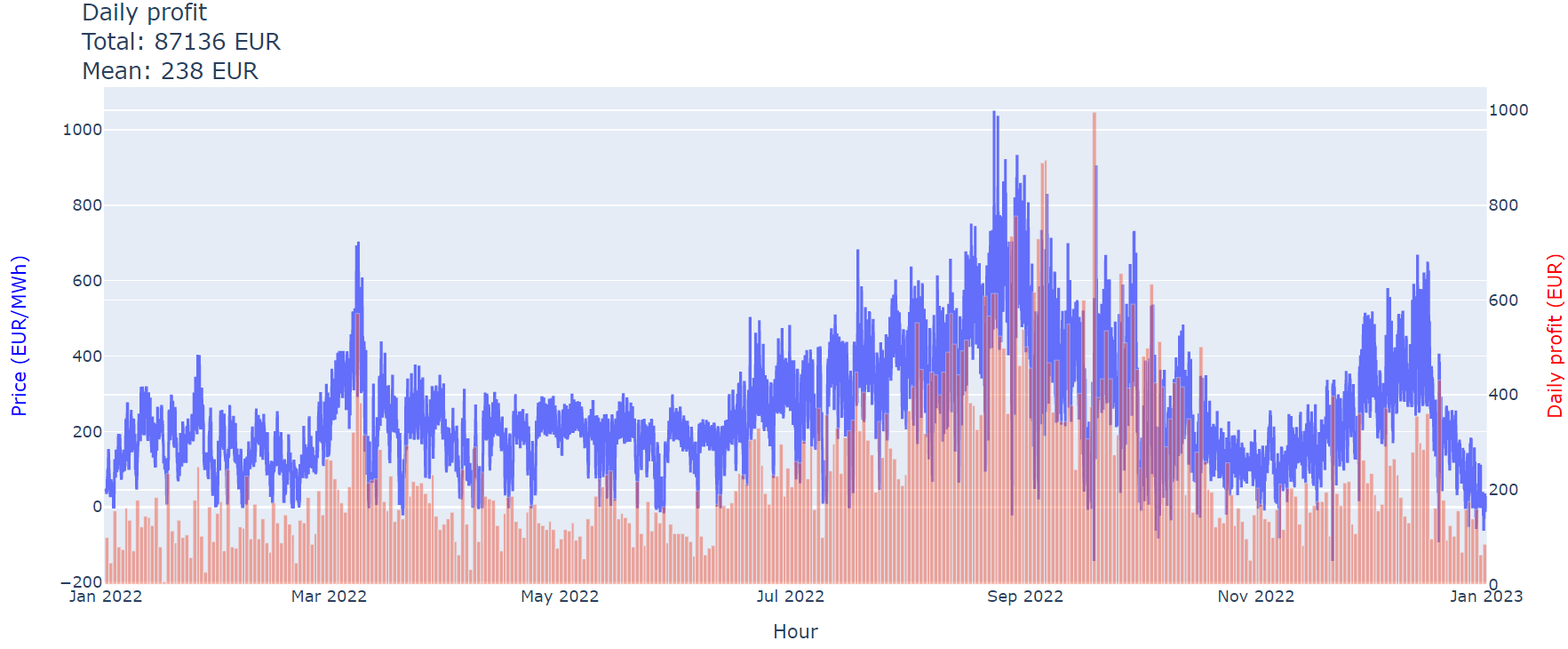}
\caption{\label{fig:full_year_pred_optim_profit} Profit distributions using MILP-P (top) with $l=28$ and MILP-O (bottom) on the German electricity market for the whole year 2022}
\end{figure}

\begin{figure}
\centering
\includegraphics[width=\columnwidth]{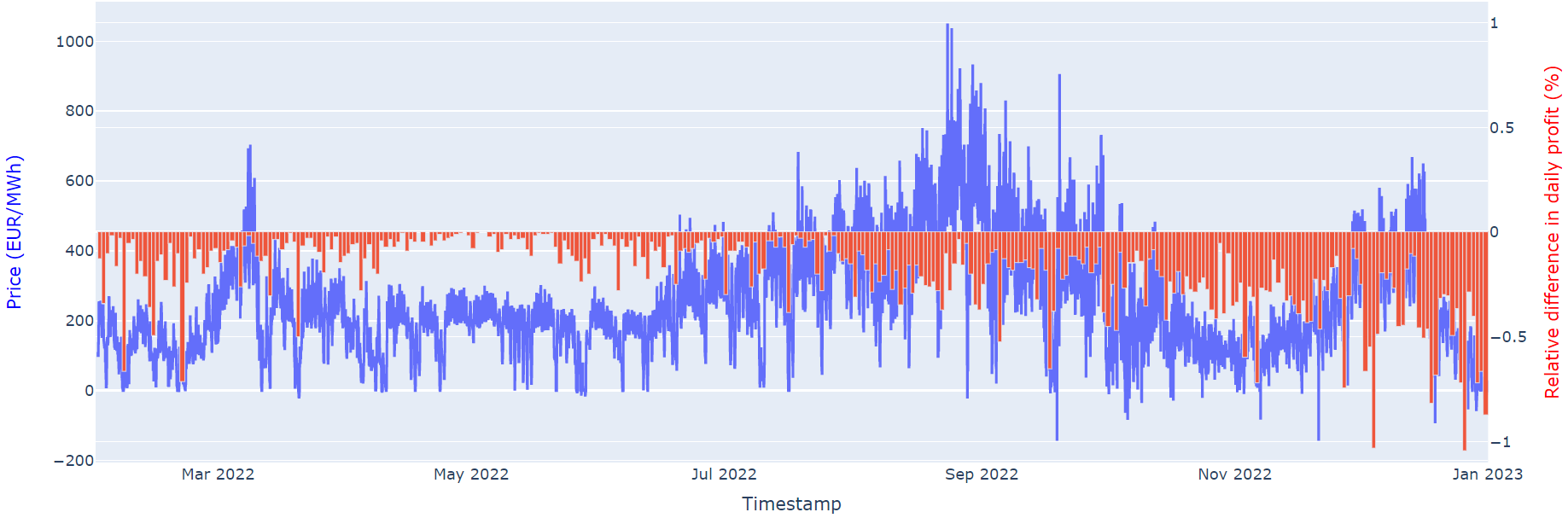}
\caption{\label{fig:difference_optim_pred} Relative difference in profit between MILP-P ($l=28$) and MILP-O (reference) on the German electricity market for the whole year 2022}
\end{figure}

The average daily profit increases with the number of days used for the computation of the price forecast $l$,  until $l=28$, for which the relative difference in profit between MILP-P and MILP-O is -19.39\%. The number of charging-discharging cycles decreases as $l$ increases, while the profit increases, meaning that the algorithm is better at tracking electricity price spreads when the forecast is based on a reasonably longer historical period. On the contrary, when $l$ is small, the schedule shows battery operations that are based on variations observed in the last $l$ days, but that do not accurately describe the usual intraday price variations. In this context, the model is overfitting on the last-day price variations. The extreme case is with $l=1$, where the model builds the schedule on the prices of the previous day, which likely showed variations that are not appearing in the current day's prices. This also increases the number of days with negative profits. We see that the mean absolute error (MAE) is not a good indicator of how good the price prediction is for the optimization task, as a low MAE is not correlated with higher profits. When $l$ is too large, e.g. $l=42$, the obtained profit decrease, as the price variation forecast becomes less specific to the current day. We show in Fig.\ref{fig:full_year_pred_optim_profit} the daily profit distribution using MILP-O and MILP-P with $l=28$. There are two main drivers of the daily profits: the daily mean electricity prices, and their intraday variations. The prices were the highest in early September, but the largest daily profits were obtained during days when the prices were lower and the intraday variations larger. Fig. \ref{fig:difference_optim_pred} shows the relative difference in profit between the two algorithms. We show the performance of both algorithms on the Danish, French, Spanish, and Italian electricity prices in Table \ref{europe} and the profits obtained with MILP-P and MILP-O on the prices for the year 2022 in other European countries in Fig. \ref{fig:profit_countries}. The smallest relative difference in profit is obtained for Italy and Denmark, where the intraday price variations vary less from one day to another.

\begin{table}
\caption{\label{table-metrics-algorithms}Results on the German electricity market for 2022}
\centering
\scalebox{0.6}{
\begin{tabular}{|c|c|c|c|c|c|c|c|}
\hline
\multirow{2}{*}{\textbf{Algorithm}} &
\multirow{2}{*}{\textbf{$l$}} &
\textbf{Avgerage}  & 
\textbf{Relative Difference in} & \multirow{2}{*}{\textbf{\#Cycles}} & \textbf{\#Days of} & 
\multirow{2}{*}{\textbf{ prediction MAE}} \\
& & \textbf{Daily Profit} & \textbf{profit with MILP-O} & & \textbf{negative profits}&
\\ \hline
\textbf{MILP-P} & 42 & 190.68 &  -20.13\%& 692.69 & 6 & 94.61 \\ \hline
\textbf{MILP-P} & 28 & 192.43 &  -19.39\%& 689.62 & 4 & 87.54 \\ \hline
\textbf{MILP-P} & 21 & 191.98 &  -19.58\%& 687.53 & 4 & 82.51 \\ \hline
\textbf{MILP-P} & 14 & 191.89 &  -19.62\%& 688.92 & 4 & 77.27 \\ \hline
\textbf{MILP-P} & 7 & 189.97 &  -20.42\%& 696.36 & 7 & 68.39 \\ \hline
\textbf{MILP-P} & 2 & 172.04 &  -27.94\%& 728.06 & 10 & 65.74 \\ \hline
\textbf{MILP-P} & 1 & 157.55 &  -34.00\%& 762.78 & 19 & 62.82 \\ \hline
\textbf{MILP-O}  &  / & 238.7 & / & 763.9  & 0 & 0 \\ \hline
\end{tabular}}
\end{table}

\begin{table}
\caption{Summary of the results obtained with MILP-P in different European countries \label{europe}}
\centering
\scalebox{0.7}{
\begin{tabular}{|c|c|c|c|c|c|c|}
\hline
\multirow{2}{*}{\textbf{Algorithm}} &
\multirow{2}{*}{\textbf{Country}} &
\multirow{2}{*}{\textbf{$l$}} &
\textbf{Avgerage}  & 
\textbf{Relative Difference in} &  
\multirow{2}{*}{\textbf{MAE prediction}} \\
& & & \textbf{Daily Profit} & \textbf{profit with MILP-O} &  
\\ \hline
\textbf{MILP-P} & Austria & 28 & 150.44 &  -28.18\%& 81.12 \\ \hline
\textbf{MILP-P} & Germany & 28 & 192.43 &  -19.39\%& 87.53  \\ \hline
\textbf{MILP-P} & France & 28 & 142.61 &  -18.48\%& 82.20 \\ \hline
\textbf{MILP-P} & Spain & 28 & 74.34 &  -16.82\%& 34.71  \\ \hline
\textbf{MILP-P} & Denmark & 28 & 157.87 &  -13.16\%& 90.75  \\ \hline
\textbf{MILP-P} & Italy & 28 & 136.70 &  -12.81\%& 68.25 \\ \hline
\end{tabular}}
\end{table}

\subsection{Impact of charging and discharging rates}
\begin{figure}
\centering
\includegraphics[width=\columnwidth]{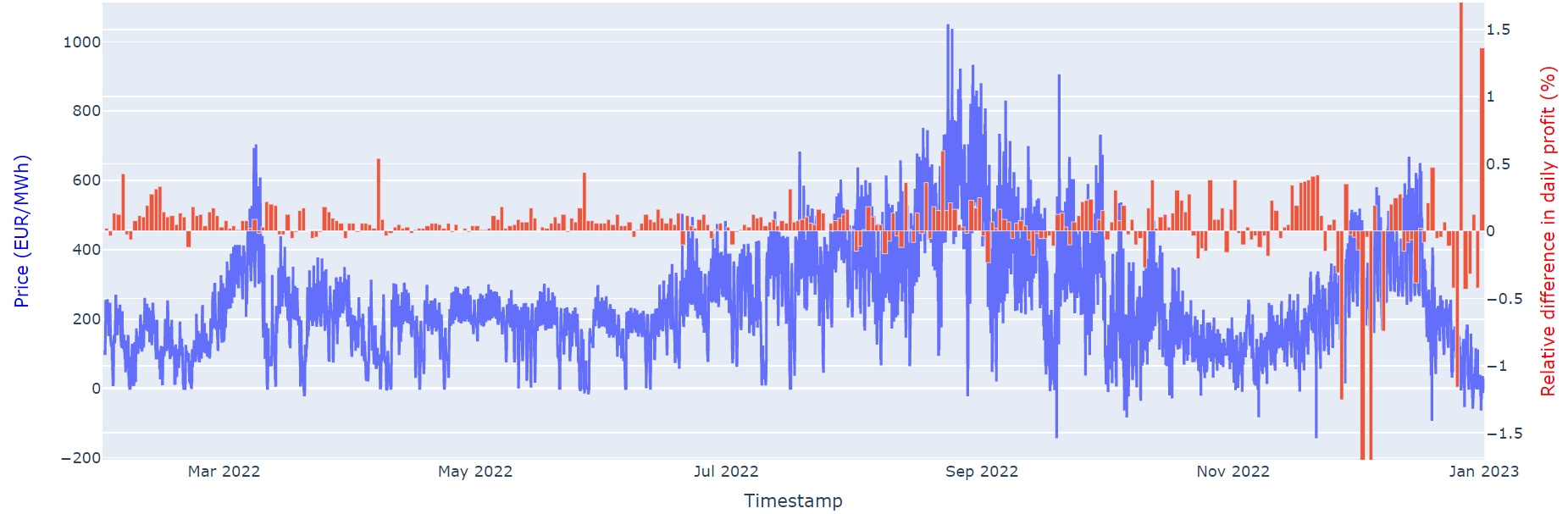}
\includegraphics[width=\columnwidth]{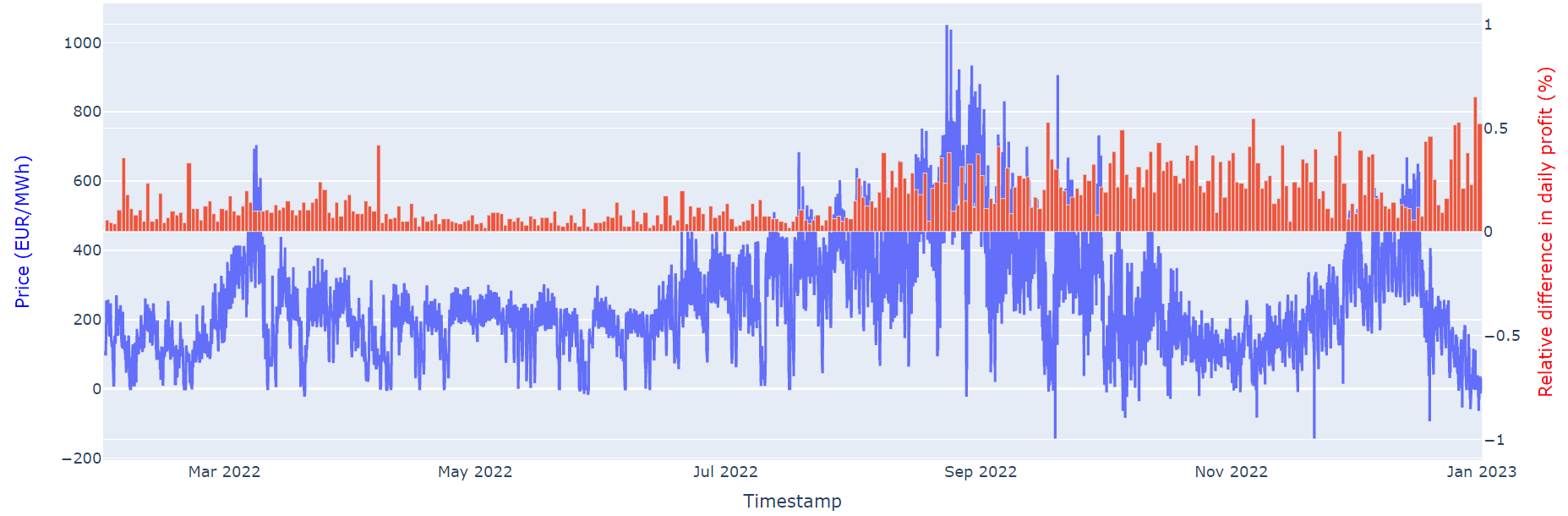}
\caption{\label{fig:relative_profit_difference} Relative difference in profit using MILP-P (top) and MILP-O (bottom) between a battery with a 1W/Wh peak charging rate and one with a 0.5W/Wh peak charging rate (reference). Both solutions are obtained using the German prices for the year 2022. Using a peak charging rate of 1W/Wh yields a profit increase of 5.86\% for MILP-P and 17.69\% for MILP-O on the entire run.}
\end{figure}

Higher charging and discharging rates increase the amount of electricity that can be bought and sold at every hour of the day. While smaller rates cause the buying and selling prices to be averaged out over multiple hours, high charging rates allow buying and selling more electricity at the local maximum and minimum prices. Fig.\ref{fig:relative_profit_difference} shows the relative difference in daily profit of MILP-P obtained with a battery with a $1$W/Wh peak charging rate and with a $0.5$W/Wh peak charging rate. The relative difference in daily profit is constant and greater than $0$ in the period when the price variability is low. In December 2022, there are days when the battery with high rates yields much lower profit than the one that has low rates. Over the entire year 2022, the difference in profit is $3806$ EUR (+5.86\% with the high-rate battery).  When the rates are high, wrong price predictions yield higher losses, since a larger volume of electricity is traded at the wrong times. This is what happens in December 2022, when the price variations are harder to predict. When the prices are known, higher rates constantly generate more profits (+17.69\% over the entire year), as shown on Fig.\ref{fig:relative_profit_difference}. Also note that assuming constant rates equal to the peak rate, as opposed to variable rates (as displayed in Fig.\ref{fig:charging_curve}), leads to an overestimation of the profits by 3.6\%.
\subsection{Impact of battery degradation and grid costs}

\begin{figure}[h]
\centering
\includegraphics[width=\columnwidth]{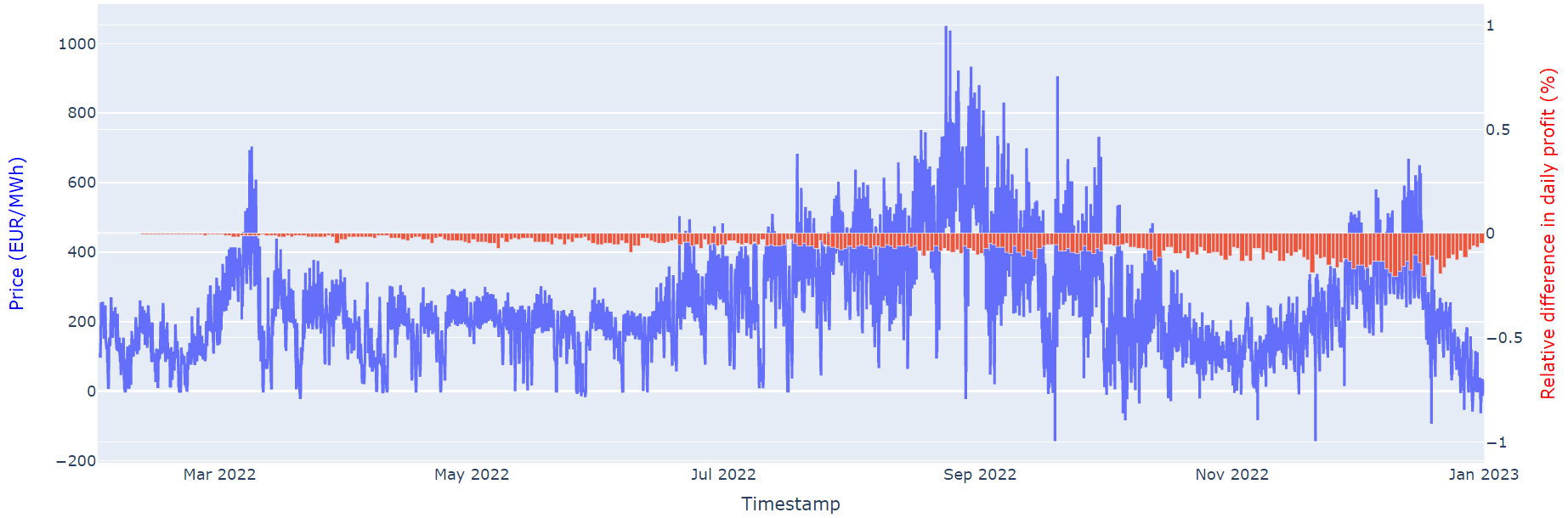}
\caption{\label{fig:battery_degradation} Relative difference in profit (6.67\% overall) due to battery degradation obtained with MILP-O}
\end{figure}

We show the relative profit delta due to capacity fading and efficiency decrease using MILP-O in Fig.\ref{fig:battery_degradation}. As the number of cycles increases, the efficiency and capacity decrease. Capacity fading reduces the amount of energy that can be stored in the battery, and a smaller efficiency decreases the amount of energy that can be taken back from the battery. We report a relative difference in profit of 6.67 \% over the entire run.

\section{Conclusion and Future Work}
\label{conclusion}
We introduced a battery model, and a Mixed Integer Linear Programming formulation of the profit maximization task of energy arbitrage on the day-ahead market. Our battery model takes into account the capacity fading effect, the decrease in discharge efficiency, and uses variable charging and discharging rates. We developed a library\footnote{The code used to run the simulations is available at https://github.com/albanpuech/Energy-trading-with-battery} to generate daily charging schedules by solving the optimization problem using forecasted prices. Its ability to take into account availability constraints, variable and fixed grid costs, custom charging curves, and battery degradation in the optimization step makes it particularly suitable for industrial use cases.
It achieved 80\% of the maximum obtainable profits on a one-year-long simulation on the 2022 German electricity prices. We report similar or better results for France (81\%) or Italy, and Denmark (87\% for both). This represents an average daily profit of 190 EUR on the German DAM\footnote{battery of 1MWh capacity, 0.5W/Wh peak charging rate, and variable grid costs of 5 EUR/MWh}. We also discussed the impact of the battery degradation on the daily profits and the profit changes when using higher charging and discharging rates. We reported a relatively small profit increase (5.86\%) when doubling the peak charging rate from 0.5W/Wh to 1 W/Wh.
In the future, we plan on improving the forecasting models. This will likely be done by using machine learning forecasting models and by formulating a loss function that better translates the relevance of the prediction for the optimization task. 

\bibliography{bibli}
\bibliographystyle{unsrt}

\end{document}